\documentclass[10pt]{article}
\usepackage{amssymb}
\newcommand{\be}{\begin{eqnarray}}
\newcommand{\ee}{\end{eqnarray}}
\newcommand{\bee}{\begin{eqnarray*}}
\newcommand{\eee}{\end{eqnarray*}}
\newcommand{\A}{{\cal A}}
\renewcommand{\d}{\mbox{d}}
\begin{document}

\begin{tabbing}
\hspace*{12cm} \= GOET-TP 66/97  \\
            \> July 1997   \\
\end{tabbing}
\vskip1.cm
%\vspace*{2cm}

\begin{center}
   \Huge \bf  Connes' distance function on one-dimensional lattices
\end{center}

\vskip1cm

\begin{center}
      {\large \bf Aristophanes Dimakis$^a$} and {\large \bf Folkert M\"uller-Hoissen$^{b,c}$} 
       \vskip.3cm 
      \begin{minipage}{14cm}
      $^a$Department of Mathematics, University of the Aegean,
      GR-83200 Karlovasi, Samos, Greece
      \vskip.1cm 
      \noindent
      $^b$Institut f\"ur Theoretische Physik,
      Bunsenstr. 9, D-37073 G\"ottingen, Germany 
      \vskip.1cm
      \noindent
      $^c$MPI f\"ur Str\"omungsforschung,
      Bunsenstr. 10, D-37073 G\"ottingen     
      \end{minipage}
\end{center} 
\vskip1cm

\begin{abstract}
\noindent
We show that there is an operator with a simple geometric significance which
yields the ordinary geometry of a linear equidistant lattice via Connes' distance
function.
\end{abstract}
\vskip1cm

\noindent
According to Connes (see \cite{Conn94} and references therein) the geodesic 
distance function
\be
     d(p,q) = \mbox{infimum of length of paths from $p$ to $q$}
\ee  
on a Riemannian manifold $M$ can be reformulated as
\be
     d(p,q) := \mbox{sup} \lbrace | f(p) - f(q) | \; ; \,
                      f  \in \A, \, \| \lbrack \hat{\cal D} ,  \hat{f} \rbrack \|
                      \leq 1 \rbrace   
\ee 
where $\A$ is a (suitably restricted) algebra of functions on $M$ represented
as multiplication operators $\hat{f}$ on a Hilbert space $\cal H$ and $\hat{\cal D}$
is the Dirac operator. The latter formulation can also be applied to discrete spaces 
and even generalized to `noncommutatice spaces'. A suitable 
replacement for the operator $\hat{\cal D}$ has to be found, however.
\vskip.2cm

In \cite{BLS94,Atzm96} a one-dimensional lattice has been considered with the
choice
\be
        (\hat{\cal D}_{s.d.} \Psi)_k = {1 \over 2 i} \, (\Psi_{k+1} - \Psi_{k-1}) 
        \qquad  k \in \mathbb{Z}    \; .         \label{BLS}
\ee
The distances calculated with this {\em symmetric difference} operator turned 
out to be given by
\be
      d(0,2n-1) = 2 \, n \, , \qquad  d(0,2n) = 2 \, \sqrt{n(n+1)}  
                                       \qquad ( n \in \mathbb{N} )  
\ee
which looks quite remote from the expected result for a linear equidistant lattice.
\vskip.2cm

In the following we show that there is another operator which 
actually produces the expected result. For simplicity, we consider a {\em finite}
set of $N$ points. Then $\A$ is the algebra of all complex functions on it.
$f \in \A$ will be represented by
\be                     
      f   \quad \mapsto \quad  \hat{f} = \left( \begin{array}{cccccc} 
          f_1 &   0        &         &       &            &    0 \\
           0   & \ddots   &          &       &            &       \\
               &               & f_N  &       &           &        \\
                &               &        & f_1  &           &        \\
                &               &        &        &  \ddots &      \\
            0  &               &        &        &            & f_N 
                                              \end{array} \right )   
\ee 
where $f_k := f(k)$ (numbering  the lattice sites by $1, \ldots, N$).
We choose the operator
\be                            \label{alternative_dirac}
        \hat{\cal D}_N :=  \left( \begin{array}{cc} 0 & {\cal D}_N^\dagger \\
                               {\cal D}_N & 0 \end{array} \right)       
\ee
on ${\cal H} = \mathbb{C}^{2N}$ where
\be
      {\cal D}_N = \left( \begin{array}{cccr} 0 & 1 & 0 & \cdots  0 \\
                                                                 \vdots & \ddots & \ddots & \vdots \\
                                                                 \vdots &  & \ddots & 1 \\
                                                                  0 & \cdots & \cdots  & 0
                                    \end{array} \right)   \; .
\ee
Then $(\A, {\cal H}, \hat{\cal D}_N)$ is a {\em spectral triple},
a basic structure in Connes' approach to noncommutative geometry 
\cite{Conn95} (see also \cite{Conn96} for a refinement).
It is called {\em even} when there is a grading operator. 
In the case under consideration such an operator is given by
\be
    \gamma := \left( \begin{array}{cc}  {\bf 1} & 0 \\ 0 & - {\bf 1}
                     \end{array} \right)  \; .
\ee
It is selfadjoint and satisfies
\be
      \gamma^2 = {\bf 1} \, , \qquad
      \gamma \, \hat{\cal D}_N = - \hat{\cal D}_N \, \gamma \, \qquad 
      \gamma \, \hat{f} = \hat{f} \, \gamma   \; .
\ee 
Let us now turn to the calculation of the distance function. With a complex
function $f$ we associate a real function $F$ via
\be
     F_1 := 0 \, , \quad F_{k+1} := F_k + | f_{k+1} - f_k | \qquad 
                k = 1, \ldots, N-1  \; .
\ee
Then $ | F_{k+1} - F_k | = | f_{k+1} - f_k | $ and
\be
      \| \lbrack \hat{\cal D}_N , \hat{f} \rbrack \, \psi \|
  = \| \lbrack \hat{\cal D}_N , \hat{F} \rbrack \, \psi \| 
\ee
for all $\psi \in \mathbb{C}^{2N}$.
Consequently, in calculating the supremum over all functions $f$ in
the definition of Connes' distance function, it is sufficient to consider only
{\em real} functions. Then $Q_N := i \, \lbrack \hat{\cal D}_N , \hat{f} \rbrack$ 
is hermitean and its norm is given by the maximal absolute value
of its eigenvalues. Instead of $Q_N$ it is simpler to consider
\be
     Q_N \, Q_N^\dagger = \mbox{diag}\left( 0, (f_2-f_1)^2, \ldots,
     (f_N-f_{N-1})^2, (f_2-f_1)^2, \ldots, (f_N-f_{N-1})^2, 0 \right)
\ee
which is already diagonal. This implies
\be
   \| \lbrack \hat{\cal D}_N , \hat{f} \rbrack \| = \mbox{max} \, \lbrace | f_2-f_1 | , 
   \ldots, | f_N-f_{N-1} | \rbrace
\ee
from which we conclude that $ d(k,l) = | k - l | $.
\vskip.2cm
                                                                          
The choice (\ref{BLS}) for $\hat{\cal D}$ was motivated by a simple discretization
procedure (which is known to cause the problem of fermion doubling in lattice
field theories). There is, however, no reason why this operator must yield the 
plain geometry of a linear equidistant lattice via Connes' distance function.
There are many geometries which can be assigned to a discrete set and these
should correspond to the choice of some operator $\hat{\cal D}$. Now it is
certainly of interest to know what distinguishes our choice (\ref{alternative_dirac}).
This is built from the operator $\cal D$ in such a way that $\hat{\cal D}$ is
selfadjoint. Moreover, the construction guarantees that there is a grading operator.
So we are left to understand the significance of $\cal D$. This matrix is the
adjacency matrix of the oriented linear lattice graph (see Fig. 1).

\begin{minipage}{5cm}
\unitlength1.cm
\begin{picture}(5,1)(-1.5,0)
\thicklines
%\linethickness{0.3mm}
%
\multiput(0,0)(1,0){6}{\circle*{0.15}}
\multiput(0,0)(1,0){3}{\vector(1,0){0.9}}
\put(3.3,0){$\ldots$}
\multiput(4,0)(1,0){1}{\vector(1,0){0.9}}
\end{picture}
\end{minipage}
\hfill
\begin{minipage}[t]{5cm}
\begin{center}
{\bf Fig. 1}
\vskip.1cm \noindent
\small
An oriented linear lattice graph. 
\end{center}
\end{minipage}    

\vskip.3cm
\noindent
This digraph plays a basic role in a formulation of lattice theories in the 
framework of noncommutative geometry \cite{DMH94_ddm} (see also
the references cited there).
\vskip.3cm
\noindent
{\em Remark.} Instead of using $\hat{\cal D}$ to define Connes' distance
function, we may use directly $\cal D$ (which, in general, is not symmetric) and no 
doubling in the representation of $f$. A simple calculation in the case treated
above shows that
\be
    \| \lbrack {\cal D}_N , f  \rbrack \| = \| \lbrack \hat{\cal D}_N , \hat{f} \rbrack \| 
\ee 
so that we obtain the same distances as before.  \hfill  {\large $\blacksquare$}
\vskip.3cm

Let us now turn to a {\em closed} linear lattice. Connecting in addition the last with 
the first point in the oriented digraph in Fig. 1, the adjacency matrix becomes
\be
      {\cal D}_{Nc} = \left( \begin{array}{cccr} 0 & 1 & 0 & \cdots  0 \\
                                                                 \vdots & \ddots & \ddots & \vdots \\
                                                                 \vdots &  & \ddots & 1 \\
                                                                  1 & 0 & \cdots  & 0
                                    \end{array} \right)   \; .
\ee
For $\Psi =(\phi,\psi) \in  \mathbb{C}^{2N}$ we find
\be
      \| \lbrack \hat{\cal D}_{Nc} , \hat{f} \rbrack \, \Psi \|^2 
  &=& \| \lbrack {\cal D}_{Nc} , f  \rbrack  \, \psi \|^2 
            + \| \lbrack {\cal D}_{Nc}^\dagger , f  \rbrack  \, \phi \|^2  \nonumber \\
  &=& \sum_{k=1}^N | f_{k+1} - f_k |^2 \, ( |\psi_{k+1}|^2 + | \phi_k|^2 ) \nonumber \\
  &\geq& \sum_{k=1}^N   | \, |f_{k+1}- a| - | f_k - a| \, |^2 
               \, ( |\psi_{k+1}|^2 + | \phi_k|^2 )  
    =  \| \lbrack \hat{\cal D}_{Nc} , \hat{F} \rbrack \, \Psi \|^2 
\ee
where $F_k = | f_k - a|$. Here and in the following, an index $N+1$ has to be replaced by $1$. 
Choosing $a = f_1$, we have $  \| \lbrack \hat{\cal D}_{Nc} , \hat{F} \rbrack \| \leq 
  \| \lbrack \hat{\cal D}_{Nc} , \hat{f} \rbrack \| $ and $F_k = | f_k - f_1 |$. It follows
that
\be 
    d(1,n) = \mbox{sup} \lbrace |F_n| \, ; \; F \mbox{ real}, \, F_1 = 0, \,  
                 \| \lbrack \hat{\cal D}_{Nc} , \hat{F} \rbrack \| \leq 1 \rbrace
\ee
The condition $ \| \lbrack \hat{\cal D}_{Nc} , \hat{F} \rbrack \| \leq 1$ is equivalent
to $ | F_{k+1} - F_k | \leq 1 $, $k=1, \ldots, N$. Let $n-1 \leq N-n+1$. It is then possible
to set the first $n-1$ terms in the identity
\be                \label{identity1}
   \underbrace{(F_2-F_1) + \cdots + (F_n - F_{n-1})}_{\mbox{$n-1$ terms}} 
 + \underbrace{(F_{n+1}-F_n) + \cdots + (F_1-F_N)}_{\mbox{$N-n+1$ terms}} = 0
\ee
each separately to $1$. Using the trivial identity
\be
     | F_n | = | (F_2 - F_1) + (F_3 - F_2) + \cdots + (F_n - F_{n-1}) |  
\ee
we now find $d(1,n) = n-1$. If $n-1 > N-n+1$, then each of the last $N-n+1$ terms in 
(\ref{identity1}) can be set to $1$. Using 
\be
  | F_n | =  | (F_1 - F_N) + (F_N - F_{N-1}) + \cdots + (F_{n+1} - F_n) |
\ee
we find $d(1,n)=N-n+1$.
\vskip.2cm

In Connes' noncommutative geometry the commutator $ \lbrack \hat{\cal D} , \hat{f} \rbrack$
represents a `differential' $\d f$. The inequality which appears in the definition of the distance
function can then be written as $\| \d f \| \leq 1$. Given a differential calculus (in the abstract
algebraic sense), in order to have a distance function we need a definition of the norm
of $\d f$. Connes defines it via a representation of the (first order) differential algebra. In
the case of discrete sets, it is natural to define a norm by
\be
    \| \d f \| = \mbox{sup} \lbrace | f(k)-f(l)| /\epsilon_{kl} \, ; \;  (kl) \in E  \rbrace
\ee
where a (di)graph structure has been assigned to the set by the first order differential calculus
(see \cite{DMH94_ddm} for details) and $E$ denotes the set of its arrows. The constants 
$\epsilon_{kl}$ are introduced here for convenience if one is interested in a continuum limit. 
The distance function is then taken to be 
\be
     d(p,q) := \mbox{sup} \lbrace | f(p) - f(q) | \; ; \,
                      f  \in \A, \, \| \d f \| \leq 1 \rbrace    \; .
\ee 
This receipe now reproduces the ordinary distances on the underlying graph.\footnote{The
distance evidently does not depend on the orientation of arrows in the digraph.}  

\vskip.2cm
\noindent
Further results on distance functions {\`a} la Connes on finite sets will be published 
elsewhere \cite{DMH97}.

\end{document}